\documentclass[oneside]{amsart}
\usepackage[latin1]{inputenc}
\usepackage[T1]{fontenc}
\usepackage{amstext}
\usepackage{amsmath, amssymb}
\usepackage[english]{babel}
\usepackage{amsfonts}

\newcommand{\R}{\mathbb{R}}
\newcommand{\C}{\mathbb{C}}
\newcommand{\A}{\mathbb{A}}
\newcommand{\N}{\mathbb{N}}
\newcommand{\K}{\mathbb{K}}

\newcommand{\E}{\mathbb{E}}

\newcommand{\Z}{\mathbb{Z}}
\newcommand{\Q}{\mathbb{Q}}
\newcommand{\pp}{\mathbb{P}}
\newcommand\kB{\mathcal{B}}

\newcommand\kR{\mathcal{R}}
\newcommand\kO{\mathcal{O}}
\newcommand\kP{\mathcal{P}}

\newcommand\kG{\mathcal{G}}
\newcommand\kH{\mathcal{H}}

\newcommand\kp{\mathfrak{p}}

\newtheorem {lem} {Lemma} [section]
\newtheorem {prop} {Proposition} [section]

\newtheorem {theo} {Theorem} [section]
\newtheorem {cor} {Corollary} [section]

\newcommand\ran{\rangle}
\newcommand\lan{\langle}

\title{The Poisson boundary of triangular matrices in a number field}
\author{Bruno Schapira}

\begin{document}

\begin{abstract} The aim of this note is to describe the Poisson
boundary of the group of invertible triangular matrices with
coefficients in a number field. It generalizes to any dimension
and to any number field a result of Brofferio \cite{Bro}
concerning the Poisson boundary of random rational affinities.
\end{abstract}

\maketitle


\bigskip

\noindent \textbf{A.M.S. classification:} 22D40; 28D05; 28D20;
60B15; 60J10; 60J50.

\section{Introduction}
The Poisson boundary is a measure space which describes the
asymptotic behavior of random walks on groups. In the same time it
gives information on the geometry of the group and provides a
representation of bounded harmonic functions (we refer for
instance to \cite{Fur} \cite{F} \cite{F2} or \cite{Kai} for more
details). Our aim in this paper is to explore the case of a group
of matrices with coefficients in a number field. More precisely we
study the group of upper diagonal matrices with non zero diagonal
coefficients.

\vspace{0.2cm} This example was treated previously by Brofferio
\cite{Bro} for matrices of size $2$ and rational coefficients,
which corresponds to the case of rational affinities. She proved
that the Poisson boundary is the product over all prime numbers
$p$ (including $p=\infty$) of "local" boundaries $C_p$, which are either a $p$-adic line, or a point.
This can be determined explicitly in function of the random walk. In
each case one can see $C_p$ as a subspace of the $p$-adic
projective line, which is the Poisson boundary of the
group of $p$-adic affinities (cf \cite{Eli} \cite{CKW}). The goal
of this paper is to generalize this result in higher dimension
$d\ge 2$. In other words we will prove that the Poisson boundary
is a product of local factors $C_p$, where for every $p$, $C_p$ is a subspace of the Furstenberg boundary, which is also the space
of flags on $\Q_p^d$. It is known that this space, or a
quotient, is the Poisson boundary of a large class of random walks
on groups of real matrices (see e.g. \cite{F} \cite{Led} or
\cite{Rau}). There is a well known decomposition of the
Furstenberg boundary called Bruhat decomposition, which coincides
for $d=2$ with the decomposition of the projective line into a line
and a point. So we will prove that each $C_p$ is a component, also
called a Bruhat cell, of this decomposition, that we determine in
function of the random walk.

Our proof follows very closely the general strategy of Brofferio \cite{Bro} in dimension $2$. However, as the technical details are a bit different, we will repeat all arguments here. So this paper can be read independently of \cite{Bro}. The main tools are the law of large numbers (for contraction) and Kaimanovich's entropy criterion (for maximality).

Such factor decompositions of the Poisson boundary were already
observed so far. For instance Bader and Shalom \cite{BS} proved
recently a general factor theorem in an adelic setting. It is in
fact rather likely that our result should extend to a more general
class of groups, such as $SL_d(\Q)$.
For this our proof "with hands" should probably be replaced by
more powerful tools, such as the Oseledec' theorem (see for
instance its use by Ledrappier \cite{Led} for the study of
discrete subgroups of semisimple groups), or a geometric argument (using for instance Kaimanovich's strip approximation criterion).

\vspace{0.2cm} \noindent \textbf{Acknowledgments:} I would like to
thank Philippe Bougerol, Sara Brofferio, Yves Guivarc'h and Vadim
Kaimanovich for enlightening discussions. I am particularly
indebted to Uri Bader for his comments on a previous version of
this paper. In particular he draw my attention to the fact that the local parts of the Poisson boundary were certainly Bruhat cells. I wish also to thank the referee for his careful reading and criticism, and for numerous advices concerning the organization of the paper.

\section{Statement of results}
\label{secresult} Let $\K$ be some number field. The reader non familiar with this may think to the particular case of the rational field $\K=\Q$ (see also section \ref{secnf} for more details). Let $\kP$ be the set of prime ideals (the set of prime numbers in the case $\K=\Q$). For $\kp\in \kP$, we denote by $|\cdot|_{\kp}$ the associated norm on $\K$ and by $\K_{\kp}$ the associated completion of $\K$. Let $\mu$ be some measure on $A(\K)$, the space of upper triangular matrices with coefficients in $\K$ and non zero diagonal coefficients. For $\kp\in \kP$ and $i\in [1,\dots,d]$, we set
\begin{eqnarray}
\label{phip}
\phi_{\kp}(i):=\int_{A(\K)} \ln |a_{i,i}|_{\kp}\ d\mu(a).
\end{eqnarray}
We denote by $B_{\kp}$ the space of flags in $\K_{\kp}^d$.  
Let $W$ be the group of permutations of the set $\{1,\dots,d\}$. 
Let $w$ be the unique element in $W$ such that $w(i)>w(j)$ if
$i<j$ and $\phi_\kp(i)\ge \phi_\kp(j)$. Let $C_\kp(\mu)$ be the 
Bruhat cell associated to $w$ in the Bruhat decomposition of $B_\kp$ (see next section for a definition). If $V_\infty$ is the set of archimedean norms on $\K$ (reduced to the usual absolute value if $\K=\Q$), then for every $v\in V_\infty$, $|\cdot|_v$, $\K_v$,..., are defined analogously.

\noindent The main result of this paper is the
\begin{theo}
\label{maintheorem} Let 
$$\mathbf{B}:=\prod_{\kp \in \kP} B_{\kp} \times \prod_{v\in V_\infty} B_v,$$ be the product of all flag manifolds. For every $\mu$ on $A(\K)$ satisfying 
\begin{eqnarray}
\label{integralcondition}
\int_{A(\K)} \sum_{i\le j}\left( \sum_{\kp\in \kP} \left|
\ln|a_{i,j}|_p\right|+\sum_{v\in V_\infty} |\ln |a_{i,j}|_v|\right)\ d\mu(a) < +\infty,
\end{eqnarray}
there exists a measure $\nu$ on $\mathbf{B}$ 
such that $(\mathbf{B},\nu)$ is the Poisson boundary of $(A(\K),\mu)$. Furthermore, $\nu$ is supported on the product of the $C_\kp(\mu)$'s and the $C_v(\mu)$'s. 
\end{theo}

We will prove this theorem in three steps, corresponding to the three following propositions:

\begin{prop}
\label{prop1}
There exists a measure $\nu$ on $\mathbf{B}$, such that the measure space $(\mathbf{B},\nu)$ is a $\mu$-boundary. 
\end{prop}

\begin{prop}
\label{prop2}
If $\mu$ satisfies \eqref{integralcondition}, then $\mu$ has finite entropy. 
\end{prop}

\begin{prop}
\label{prop3}
For $\nu$-almost all $z\in \mathbf{B}$, the asymptotic entropy $h^z$ of the conditional measure $\pp^z$ vanishes.
\end{prop}

\noindent We will recall all necessary definitions about entropy in the next section. Assuming these propositions, Theorem \ref{maintheorem} is then a consequence of Kaimanovich's criterion: 

\begin{theo}[Kaimanovich \cite{Kai} Theorem 4.6]
Let $G$ be a countable group and $\mu$ a probability measure on $G$ with finite entropy. Then a $\mu$-boundary $(\mathbf{B},\nu)$ is the Poisson boundary if, and only if, for $\nu$-almost all $z\in \mathbf{B}$, the asymptotic entropy $h^z$ of the conditional measure $\pp^z$ vanishes.  
\end{theo}

\vspace{0.2cm} Let us describe now the organization of the paper.
In the next section we detail our notations and recall some
preliminary background on $\mu$-boundaries, Poisson boundary, entropy and Bruhat decomposition. Then we prove the three propositions above in the particular case of the rational field, which is easier in a first reading. 
Section \ref{secmuboundary} is devoted to the proof of Proposition \ref{prop1}, section \ref{secgauge} to Proposition \ref{prop2}, and section \ref{secestimate} to Proposition \ref{prop3}. The
last section is devoted to the case of number fields. There are
some adjustments to make in the proof that we explain. Finally the appendix is devoted to the proof of a technical result.

\section{Preliminaries}
Let $d\ge 1$ be an integer. Recall that $A(\K)$ is the set of upper
triangular matrices of size $d$ with coefficients in $\K$ and non
zero diagonal coefficients. So if $a=(a_{i,j})_{i,j} \in A(\K)$,
we have $a_{i,j}=0$, if $i>j$, and $a_{i,i}\neq 0$ for $1\le i\le
d$. For $n\ge 1$ and
$z=(z_1,\cdots,z_n)\in \K_\kp^n$, we set
$|z|_\kp=\max_{i=1,\cdots,n}|z_i|_\kp$. 

\vspace{0.2cm} \noindent \textbf{Random walk, $\mu$-boundaries, and Poisson
boundary:} Let $\mu$ be a measure on $A(\K)$. We consider a
sequence $(g_n)_{n\ge 1}$ of i.i.d. random variables of law $\mu$
on $A(\K)$. The random walk $(x_n)_{n\ge 0}$ of law $\mu$ on
$A(\K)$ is defined by
$$x_n:=g_1\cdots g_n.$$
We denote by $\pp$ the law of $(x_n)_{n\ge 1}$ on the path space
$A(\K)^\N$.

Assume that $B$ is a locally compact space, endowed with a
measure $\nu$ and an action of $A(\Q)$. We say that $\nu$ is
$\mu$-invariant (also known as $\mu$-stationary or $\mu$-harmonic), if
$$\int_{A(\K)} (g\nu)d\mu(g)=\nu,$$ where for
all $g\in A(\K)$, $g\nu$ is defined by
$$g\nu(f)=\int_{\mathbf{B}}f(gz) d\nu(z),$$ for all continuous
functions $f$. In this case, according to Furstenberg \cite{F,F2}, we
say that $(B,\nu)$ is a $\mu$-boundary if, $\pp$-almost
surely $x_n\nu$ converges vaguely to a Dirac measure. 
Then the Poisson boundary $(\mathbf{B},\nu)$ is defined as the maximal $\mu$-boundary, 
i.e. it is a $\mu$-boundary such that any other $\mu$-boundary is one of its measurable $G$-equinvariant quotients.  
For any element $\mathbf{x}=(x_n)_{n\ge 1} \in A(\K)^\N$, we define $\mathbf{bnd}\ \mathbf{x}\in \mathbf{B}$ by
$$\lim_{n\to +\infty} x_n \nu = \delta_{\mathbf{bnd}\ \mathbf{x}}.$$
If $z\in \mathbf{B}$, it is possible to define (in the sense of Doob, see \cite{Kai} for details), the law $\pp^z$ of $(x_n)_{n\ge 1}$ conditioned by $\mathbf{bnd}\ \mathbf{x}=z$. Then for $n\ge 0$, $\pp_n^z$ denotes the projection of $\pp^z$ on the $n^{\textrm{th}}$ coordinate. In fact if $(B,\nu)$ is any $\mu$-boundary and $z\in B$, $\pp^z$ and $\pp^z_n$ are defined similarly.

Let $\text{sgr}(\mu)$ be the semi-group generated by the support
of $\mu$, i.e. $\text{sgr}(\mu)=\cup_n \text{supp}(\mu^{*n})$. We
say that a function $f$ on $\text{sgr}(\mu)$ is $\mu$-harmonic if
$$\int_{A(\K)} f(gg')d\mu(g')=f(g),$$ for all $g\in
\text{sgr}(\mu)$. Furstenberg \cite{F} proved that there is an
isometry between the space $\kH^\infty(A(\K),\mu)$ of bounded
$\mu$-harmonic functions $f$ on $\text{sgr}(\mu)$ and the space
$L^\infty(\mathbf{B},\nu)$ of bounded functions $F$ on
$\mathbf{B}$. The isometry is given by the formula
$$F(\mathbf{bnd}\ \mathbf{x})=\lim_{n\to\infty} f(x_n),
\qquad f(g)=\int_{\mathbf{B}} F(g z)\ d\nu(z).$$ The second
formula is the so called Poisson integral representation formula
of bounded harmonic functions.

\vspace{0.2cm} \noindent \textbf{Entropy and asymptotic entropy:} 
The entropy of a measure $\mu$ on a countable group $G$ is given by the formula: 
$$H(\mu)=-\sum_{g\in G} \mu(g)\ln \mu(G).$$ 
If $(B,\nu)$ is a $\mu$-boundary and $z\in B$  the conditional asymptotic entropy $h^z$ is defined by: 
$$h^z:=-\lim_{n\to +\infty} \frac{\log \pp_n^z(x_n)}{n}.$$

\vspace{0.2cm} \noindent \textbf{Some structure and the Bruhat
decomposition:} Let $G=GL_d$. We denote by $\Delta$ the set of
diagonal matrices with non zero diagonal coefficients. We denote
by $\delta=\textrm{diag}(\delta_1,\dots,\delta_d)$ the diagonal
matrix with entries $\delta_{i,i}=\delta_i$, $i=1,\dots,d$. Let
$U$ be the group of upper triangular matrices with one's in the
diagonal (unipotent matrices). The notation $U(\kR)$, where $\kR$
is some ring, means that the coefficients strictly upper the
diagonal are in $\kR$. Let $\overline{U}$ be the group of lower
triangular matrices with one's in the diagonal. We set $A=\Delta
U$ and $\overline{A}=\Delta \overline{U}$. We denote by $W$ the
Weyl group, identified with the subgroup of permutation matrices. Its
action by conjugation on $\Delta$ permutes the coordinates of the
diagonal. In this way $W$ can also be identified with the group
of permutations of the set $\{1,\dots,d\}$. For $w\in W$, we set
$U_w=w\overline{U}w^{-1} \cap U$ and $U^w=wUw^{-1}\cap U$. We have
$U=U^wU_w$ and $U^w\cap U_w= \{\textrm{Id}\}$. An element $u\in U$ lies in $U^w$ if, and
only if, $u_{i,j}=0$ whenever $i<j$ and $w(i)>w(j)$.

The \textit{Bruhat decomposition} (see e.g. \cite{GR} or
\cite{War}) says that $G$ can be decomposed in the following
disjoint union:
$$G=\coprod_{w\in W} Aw\overline{A}=\coprod_{w\in W} U^ww\overline{A}.$$
The components $Aw\overline{A}$ are called the Bruhat cells. In
the quotient space $G/\overline{A}$ they are identified with the
groups $U^w$ by the map 
$$U^w \to G/\overline{A}, \quad u\mapsto uw\overline{A}.$$

\noindent For $\kp \in \kP$, resp. $v\in V_\infty$, we recall that $C_\kp(\mu)$, resp. $C_v(\mu)$, denotes the cell $U^w(\K_\kp)$, resp. $U^w(\K_v)$, associated to the $w\in W$ such that $w(i)>w(j)$ if $i<j$ and $\phi_\kp(i)\ge \phi_\kp(j)$, resp. $\phi_v(i)\ge \phi_v(j)$. If $u\in U(\K_\kp)$, resp. $u\in U(\K_v)$, we denote by $\overline{u}$ its component in this $U^w(\K_\kp)$, resp. $U^w(\K_v)$, according to the product decomposition $U=U^wU_w$.

\vspace{0.1cm} \noindent The action of $a\in A(\K)$ on $b\in C_p(\mu)$
is defined as follows. Assume that $a=u\delta$ with $u\in U$ and
$\delta \in \Delta$. Then $$a\cdot b := \overline{a b
\delta^{-1}}.$$ In other words $a\cdot b$ is the unique element of
$U^w$ representing $ab$ in $G/\overline{A}$.

\section{Proof of Proposition \ref{prop1}} \label{secmuboundary} To simplify a little the notations and the arguments, we assume here and in the next three sections, that $\K=\Q$. In this case we denote by $\kP^*$ the union of $\kP$, the set of prime numbers, and $\{\infty\}$, which corresponds to the usual absolute value. So the notation $\Q_\infty$ denotes the field of reals $\R$, and $|\cdot|_\infty$ the absolute value on $\R$.  
Let $p\in \kP^*$. Let $(e_1,\dots,e_d)$ be the canonical
basis of $\Q_p^d$. Let
$$J_d=\{i\le d\mid \phi_p(i)\ge \phi_p(d)\},$$ and let $r$ be the
cardinality of $J_d$. Assume that $j_1<\cdots<j_r=d$ are the elements
of $J_d$. We denote by $\Lambda^r_{\textrm{sub}}\Q_p^d$ the subspace of
$\Lambda^r \Q_p^d$ generated by the elements $e_{i_1}\wedge \ldots \wedge
e_{i_r}$ such that $i_s\le j_s$ for all $s\in [1,\dots,r]$. We
denote by $\kB_r$ the basis of $\Lambda^r_{\textrm{sub}}\Q_p^d$ made up
of these elements ranged in lexicographical order. Let also $m$
be the dimension of this subspace. Each $a\in A(\Q)$ defines an
endomorphism of $\Lambda^r\Q_p^d$, by setting $a(v_1\wedge \ldots \wedge
v_r):=av_1\wedge \ldots \wedge av_r$. We denote by $a^{(r)}$ the
restriction of this endomorphism to $\Lambda^r_{\textrm{sub}}\Q_p^d$.
Observe that it has a triangular matrix representation in the
basis $\kB_r$. This provides a representation of $A(\Q)$ on the
subspace $P\Q_p^m$ of $\Q_p^m$ whose vectors have last coordinate
equal to $1$. Indeed for $u\in P\Q_p^m$ and $a\in A(\Q)$, we
define the (projective) action $a\cdot u$ of $a$ on $u$ by
$$a \cdot u :=\frac{1}{\prod_{j\in J}a_{j,j}}a^{(r)}u.$$
\begin{lem}
\label{convergence} For any $u\in P\Q_p^m$, the sequence
$(x_n\cdot u)_{n\ge 1}$ converges $\pp$-a.s. Moreover the limit,
that we denote by $(Z^p_k(d))_{k\le m}$, does not depend on the
choice of $u$.
\end{lem}
\begin{proof} We assume that $m\ge 2$, otherwise there
is nothing to prove. Let $a\in A(\Q)$. We have observed that
$a^{(r)}$ has a triangular matrix representation in the basis
$\kB_r$. We put $a'=\frac{1}{\prod_{j\in J} a_{j,j}}a^{(r)}$. Then
$a'_{m,m}=1$, and for $k<m$, there exists a subset $K\subset
\{1,\dots,d\}$ of cardinality $r$, such that
$a'_{k,k}=\frac{\prod_{j\in K}a_{j,j}}{\prod_{j\in J}a_{j,j}}.$
Therefore \begin{eqnarray} \phi'_p(k)&: =& \E[\ln |a'_{k,k}|_p]=
\sum_{j\in K} \E[\ln |a_{j,j}|_p]-\sum_{j\in
J}\E[\ln |a_{j,j}|_p] \\
 &= & \nonumber \sum_{j\in K\cap J^c}
\underbrace{(\phi_p(j)-\phi_p(d))}_{< 0}-\sum_{j\in J\cap
K^c}\underbrace{(\phi_p(j)-\phi_p(d))}_{\ge 0} < 0.
\end{eqnarray}
For $n\ge 1$, let
$$x'_n=g'_1\cdots g'_n.$$ If $k<m$, by the law of large
number, a.s.
\begin{eqnarray}
\label{xprimk}
\frac{\ln|(x'_n)_{k,k}|_p}{n}=\frac{1}{n}\sum_{l=1}^n
\ln|(g'_l)_{k,k}|_p  \to \phi'_p(k),\text{  when }n\to \infty.
\end{eqnarray}
\textit{first step:} Let $k$ be fixed. Let us prove by induction
on $k'\in[k,\dots,m-1]$, that for all $k\le k'\le m-1$, there
exists a.s. $\alpha>0$ such that $\frac{1}{n}\ln
|(x'_n)_{k,k'}|_p\le -\alpha$ for $n$ large enough. If $k'=k$ the
result is immediate from (\ref{xprimk}). We assume now the result
for $l\le k'-1$ (with $k'>k$), and we prove it for $k'$. By the
induction hypothesis there exists a.s. $\alpha
>0$ and $N_1\ge 1$ such that for all $l\in[k,\dots,k'-1]$,
$$n\ge N_1 \Rightarrow |(x'_n)_{k,l}|_p\le e^{-n\alpha}.$$
Let $\epsilon>0$ be such that $\alpha -4\epsilon >0$, and
$\phi'_p(k')+\epsilon<0$. By (\ref{xprimk}) there exists a.s.
$N_2\ge N_1$ such that
\begin{eqnarray} \label{xprimnkprim} n\ge N_2\Rightarrow
e^{n(\phi'_p(k) - \epsilon)} \le |(x'_n)_{k',k'}|_p\le
e^{n(\phi'_p(k)+ \epsilon)}.
\end{eqnarray} We set $c_n=\max_{i,j}|(g'_{n+1})_{i,j}|_p$. We
have a.s. $\frac{\ln c_n}{n}\to 0$. Thus there is a.s. some
integer $N_3\ge N_2$ such that $n\ge N_3 \Rightarrow c_n\le
e^{n\epsilon}$. Finally we set
$u_n:=\frac{(x'_n)_{k,k'}}{(x'_n)_{k',k'}}$. We have for all $n\ge
1$,
\begin{eqnarray}
\label{unplus1}
u_{n+1}=u_n+\underbrace{\frac{1}{(x'_{n+1})_{k',k'}}\sum_{l=k}^{k'-1}(x'_n)_{k,l}(g'_{n+1})_{l,k'}}_{:=r_n}.
\end{eqnarray}
With the previous notations we have a.s. for $n\ge N_3$,
\begin{eqnarray}
\label{rnp} |r_n|_p\le de^{-n\phi'_p(k') -n(\alpha -3\epsilon)}.
\end{eqnarray}
We set $C=\max_{n\le N_3} |r_n|_p$. Hence by (\ref{xprimnkprim}),
(\ref{unplus1}) and (\ref{rnp}), we have a.s. for any $N\ge N_3$,
$$|(x'_N)_{k,k'}|_p\le |(x'_N)_{k',k'}|_p\sum_{n=1}^N|r_n|_p\le Nd(e^{-N(\alpha -4\epsilon)}+Ce^{N(\phi'_p(k')+\epsilon)}),$$
and the result for $k'$ follows.

\noindent \textit{second step:} We have for $n\ge 1$, and
any $k<m$,
$$(x'_{n+1})_{k,m}=(x'_n)_{k,m}+\sum_{l=k}^{m-1}(x'_n)_{k,l}(g'_{n+1})_{l,m}.$$
As a consequence $(x'_n)_{k,m}$ is the partial sum of a series
whose general term converges a.s. to $0$ exponentially fast (by
the first step). Thus it is almost surely convergent. Now take $u\in P\Q_p^m$. By definition $x_n\cdot u= x'_nu$ for all $n$. So we see that $x_n \cdot u$ converges a.s. to some $(Z_k^p(d))_{k\le m} \in P\Q_p^m$, where for every $k<m$, $Z_k^p(d)$ is the limit of $(x'_n)_{k,m}$, which is independent of $u$. This finishes the proof of the lemma.
\end{proof}

\noindent This lemma says that $P\Q_p^m$ equipped with the
law of $(Z^p_k(d))_{k\le m}$ is a $\mu$-boundary. But it implies in fact the 
\begin{cor} 
\label{localpart}
For every $p\in \kP^*$, there exists a measure $\nu_p$ on $C_p(\mu)$ such that the measure space $(C_p(\mu),\nu_p)$ is a $\mu$-boundary. 
\end{cor}
\begin{proof} For all $d'\le d$, we define the minor of size $d'$
of an element $a\in A(\Q)$ as the matrix $d'\times d'$ in the
upper left corner of $a$. These matrices act in the same way on
$\Lambda_{\textrm{sub}}^{r(d')}\Q_p^{d'}$, where $r(d')$ is the
cardinality of $J_{d'}$. Hence Lemma \ref{convergence} holds as well
in this setting. This provides new $\mu$-boundaries and new
vectors $(Z^p_k(d'))_{k\le m(d')}$, where $m(d')$ is the dimension
of $\Lambda_{\textrm{sub}}^{r(d')}\Q_p^{d'}$. We claim that the
set of vectors $(Z^p(2),\cdots,Z^p(d))$ is associated to an
element $(Z_{i,j}^p)_{1\le i,j\le d}$ of $C_p(\mu)$. More precisely we
claim that we can define the columns $Z^p_j$ of $Z^p$, for
$j=1,\cdots,d$, recursively by
$$Z^p_{i_1}\wedge \cdots \wedge Z^p_j=Z^p(j),$$
where $i_1,\cdots,j$ are the elements of $J_j$. Indeed the set,
let say $S$, of vectors $(V(2),\cdots,V(d))$ which are associated
to an element of $C_p(\mu)$ by this way is stable under the action of
$A(\Q)$. But since the limit $(Z^p(2),\cdots,Z^p(d))$ is
independent of the starting point, which can be chosen in $S$, it
must be also in $S$. Thus if $\nu_p$ is the law of the associated
$Z^p$, we get that $(C_p(\mu),\nu_p)$ is a $\mu$-boundary. 
\end{proof}

\noindent \textbf{Proof of Proposition \ref{prop1}:} It suffices to observe the elementary fact that a product of $\mu$-boundaries is a $\mu$-boundary. So if we define $\nu$ on $\mathbf{B}$ to be the law of $(Z^p)_{p\in \kP^*}$, we get from Corollary \ref{localpart} that $(\mathbf{B},\nu)$ is a $\mu$-boundary. \hfill $\square$

\section{Gauges on $A(\Q)$ and proof of Proposition \ref{prop2}}
\label{secgauge} We denote by $\A$ the adele ring of $\Q$, i.e.
the restricted product $\Pi'_{p\in \kP^*}\Q_p$ (see e.g. \cite{Wei}).
The notation $\Pi'$ means that if $(z^p)_{p\in \kP}\in \A$, then
for all $p$ but a finite number, $|z^p|_p \le 1$. Let $H$ be the
group of upper triangular matrices with non zero rational diagonal
coefficients and strictly upper diagonal coefficients in $\A$. In
other words $$H:= U(\A)\ \Delta(\Q).$$ We have a natural injection
$i_{\A}$ from $\Q$ into $\A$ and therefore also an injection
$i_{H}$ from $A(\Q)$ into $H$. Via $i_H$ we will sometimes
identify elements in $A(\Q)$ with their image in $H$. For $q\in
\Q^*$, we set
$$\lan q \ran:=\sum_{p\in \kP}|\ln|q|_p|.$$
In particular for every irreducible fraction $q=\pm r/s$ of
integers, one has $\lan q\ran=\ln r+\ln s$. If
$\delta=\textrm{diag}(\delta_1,\dots,\delta_d) \in \Delta(\Q)$, we
set
$$\lan \delta \ran:=\sum_{i=1}^d \lan \delta_i\ran.$$
For $b=(b^p)_{p\in \kP^*}\in \A$ and $p\in \kP^*$, we set
$$\lan b\ran_p^+:= \ln^+|b^p|_p,$$
where $\ln^+$ denotes the positive part of the function $\ln$ and
$$\lan b\ran^+:=\sum_{p\in \kP^*} \lan b\ran_p^+$$
If $u\in U(\A)$ and $p\in \kP^*$, we set
$$\lan u\ran^+_p:= \sum_{i<j} \lan u_{i,j}\ran_p^+ \quad \text{and} \quad \lan u\ran^+:= \sum_{i<j} \lan u_{i,j}\ran^+.$$
Let $h\in H$ and let $h=u\delta$ be its decomposition in
$U(\A)\Delta(\Q)$. We define the adelic length of $h$ by
$$||h||:=\lan u\ran^+ + \lan \delta \ran.$$ The adelic length is
not sub-additive but it is almost
the case. Indeed for any $q,q'\in \Q^*$,
$$\lan qq'\ran \le \lan q\ran +\lan q'\ran,$$ and for any $b,b',b''\in \A$, and $q\in \Q^*$,
$$\lan b+qb'b''\ran^+ \le \ln 2 + \lan b\ran^++\lan b'\ran^+ +\lan b''\ran^++ \lan q\ran.$$
Using these relations we can find constants $K>0$ and $K'>0$ such
that for all $h,h'\in H$,
$$||hh'||\le K +
K'(||h||+||h'||).$$ Now we consider the family of gauges
$(\kG_k^h)_{k\in \N}$ on $A(\Q)$ defined for $k\ge 0$ and $h\in
H$, by
\begin{eqnarray}
\label{defgauge}
\kG_k^h:=\{a\in A(\Q) \mid ||a^{-1}h||\le k\}.
\end{eqnarray}

We have 

\begin{lem}
\label{lemgauges}
The family of gauges $\{\kG^h\}_{h\in H}$ has uniform exponential
growth, i.e. there exists $C'>0$ such that $\text{Card}
\{\kG_k^h\}\le e^{C'k}$ for all $h\in H$ and all $k\in \N-\{0\}$.
\end{lem}

\begin{proof} First, since the inverse map is a
bijection of $A(\Q)$, we can always replace $a^{-1}$ by $a$ in the
definition of the gauges. Now let $h_0$ be the unit element of
$H$, and let $a\in A(\Q)$ be such that $||ah_0|| \le k$. In this
case $\lan a_{i,i} \ran\le k$ for any $i\le d$, and $\lan a_{i,j}/a_{j,j}\ran^+
\le k$ for any $i<j$. But the number of rational $q\neq 0$ such
that $\lan q\ran \le k$ is lower than $2e^{2k}$. Moreover, for any
rational $q$, $\lan q\ran^+ \ge \lan q\ran /2$. Thus
$$\lan \frac{a_{i,j}}{a_{j,j}}\ran^+ \ge \frac{1}{2}
\lan \frac{a_{i,j}}{a_{j,j}}\ran \ge \frac{1}{2}(
\lan a_{i,j}\ran-\lan a_{j,j}\ran),$$ which implies
$$\text{Card}(\kG_k^{h_0})\le (2e^{6k})^{d^2}.$$ Now let $h\in H$.
Since the multiplication by any element is a bijection of $A(\Q)$,
we do not change the cardinality of the $\kG_k^h$ if we multiply
to the left $h$ by an element in $A(\Q)$. Hence, multiplying them
if necessary by $\textrm{diag}(h_1^{-1},\dots,h_d^{-1})$ we can
always suppose that $h_1=\dots=h_d=1$. Then it is elementary to find $b\in A(\Q)$ such that $||h^{-1}b||=0$. Hence for any
$a\in A(\Q)$,
$$||ab||= ||a h h^{-1}b||\le
K+K'(||ah||+||h^{-1}b||).$$ Thus $\kG^h_k\subseteq
\kG^{b}_{K+K'k}$, which has the same cardinality as
$\kG^{h_0}_{K+K'k}$, since $b\in A(\Q)$. This concludes the proof
of the lemma. 
\end{proof}

\noindent If $\kG=(\kG_n)_{n\ge 1}$ is a gauge on a countable group $G$, and if $g\in G$, we set 
$$|g|_\kG:=\inf\{n\mid g\in \kG_n\}.$$ 
Then if $\mu$ is a measure on $G$, the first moment $|\mu|_\kG$ of $\mu$ with respect to $\kG$ is defined by: 
$$|\mu|_\kG:= \sum_{g\in G} |g|_\kG\mu(g).$$

The proof of Proposition \ref{prop2} follows now from Derriennic's criterion:  

\begin{theo}[Derriennic \cite{Der}] 
Let $\mu$ be a probability measure on a countable group $G$. If $\mu$ has finite first moment with respect to some exponentially growing gauge, then $\mu$ has finite entropy.  
\end{theo} 

Indeed Hypothesis \eqref{integralcondition} says exactly that $\mu$ has finite first moment with respect to the gauge $(\kG_n^{h_0})_{n\ge 1}$, where $h_0$ is the identity matrix, and Lemma \ref{lemgauges} assures in particular that $\kG^{h_0}$ has exponential growth.

\section{Proof of Proposition \ref{prop3}} \label{secestimate} We start by some preliminary estimates. Remember that if $q\in \Q^*$, then 
$$\lan q\ran = \sum_{p\in\kP} |\ln |q|_p|.$$
Remember also the definition of $\phi_p$ from \eqref{phip}. We have
\begin{lem}
\label{qni} For $i\le d$, and $n\ge 1$, let $q_n^i=\prod_{p\in
\kP}p^{-\left[ n\frac{\phi_p(i)}{\ln p}\right] },$ where for
$x\in \R$, $[x]$ denotes the integer part of $x$ if $x\ge 0$, and
the opposite of the integer part of $-x$ otherwise. Then
$$\frac{\lan (x_n)^{-1}_{i,i}q_n^i\ran}{n} \to 0, \text{ in } L^1.$$
\end{lem}
\begin{proof} For $p\in \kP$, the ergodic theorem implies that
$$\frac{\ln |(x_n)^{-1}_{i,i}q_n^i|_p}{n}=\frac{-\sum_{k=1}^n
\ln|(g_k)_{i,i}|_p + \ln p \left[n\frac{\phi_p}{\ln p} \right]}{n}
\to 0$$ in $L^1$. Thus by the dominated convergence theorem, the
sequence
$$\E \left[ \frac{\lan(x_n)^{-1}_{i,i}q_n^i\ran}{n} \right] = \sum_{p\in \kP} \frac{\E\left[|\ln
|(x_n)^{-1}_{i,i}q_n^i|_p|\right]}{n}$$ converges to zero because
each term of the infinite sum converges to zero and is dominated
by $\E\left[|\ln |a_{i,i}|_p|\right]+|\phi_p|$ whose series is
convergent by \eqref{integralcondition}. 
\end{proof}

\noindent Let now $P$ be some finite subset of $\kP^*$ and let
$q_n=\textrm{diag}(q_n^1,\dots,q_n^d)$. We set
\begin{eqnarray*}
\begin{array}{cccl}
\pi_n^P : & \prod_{p\in P^*} B_p & \longrightarrow & H \\
          & (z^p)_{p\in P^*} & \mapsto &
          \mathbf{z}\ q_n,
          \end{array}
\end{eqnarray*}
where for $i<j$, $\mathbf{z_{i,j}} \in \A$ is defined by
\begin{eqnarray*}
\mathbf{z^p_{i,j}}= \left\{ \begin{array}{cc}
z_{i,j}^p & \text{if }p\in P, \\
0          & \text{otherwise.}
          \end{array}
          \right.
\end{eqnarray*}
We set also $Z^P:=(Z^p)_{p\in P} \in \prod_{p\in P} B_p$. The main
result of this section is the following proposition.
\begin{lem}
\label{estimgauge} Let $P$ be some finite subset of $\kP^*$ containing $\infty$. For $p\in \kP^*$, let $K_p=\sum_{r\le
s}\int_{A(\Q)} |\ln|a_{r,s}|_p|\ d\mu(a)$. Let $\epsilon >0$ be some constant. Then there exists a
constant $C>0$, such that
$$\pp \left[\frac{||x_n^{-1}\pi^P_n(Z^P)||}{n}\le 
\epsilon+C\sum_{p\notin P}K_p \right] \to 1.$$
\end{lem}
\begin{proof} Assume that
$$x_n=u_n \delta_n \quad \forall n\ge 1,$$
with $u_n \in U$ and $\delta_n \in \Delta$. We have
$$||x_n^{-1}\pi^P_n(Z^P)||=\lan \delta_n^{-1}q_n\ran +
\lan x_n^{-1}\pi_n^P(Z^P)q_n^{-1} \delta_n\ran^+.$$ First we know by
Lemma \ref{qni} that $\lan \delta_n^{-1}q_n\ran /n$ converges to $0$ in
$L^1$. So it converges also to $0$ in probability. Next
$$\lan x_n^{-1}\pi_n^P(Z^P)q_n^{-1} \delta_n\ran^+ =\sum_{p\in P}\lan x_n^{-1}Z^p\delta_n \ran^+_p+\sum_{p\notin P}
\lan x_n^{-1}\delta_n\ran^+_p.$$

\vspace{0.2cm} \noindent \textit{first step: the sum over $p\notin
P$. }
\newline
For $i\le j$ and $N\ge 1$ we have
$$\frac{(x_{N+1}^{-1})_{i,j}}{(x_{N+1}^{-1})_{i,i}}=\frac{(x_N^{-1})_{i,j}}{(x_N^{-1})_{i,i}}+\underbrace{\sum_{k=i+1}^j\frac{(x_N^{-1})_{k,j}}{(x_{N+1}^{-1})_{i,i}}
(g_{N+1}^{-1})_{i,k}}_{:=r_N}.$$ By the ultra-metric property we
get \begin{eqnarray} \label{xNrN} \ln^+|(x_{N+1}^{-1})_{i,j}|_p\le
|\ln|(x_{N+1}^{-1})_{i,i}|_p|+\max_{1\le n\le N}\ln^+|r_n|_p.
\end{eqnarray}
For $n\ge 1$, we set $c_n=\max_{r,s}|(g_{n+1}^{-1})_{r,s}|_p$.
Observe that for some constant $C_1>0$, $E[|\ln c_n|] \le C_1K_p$.
By \eqref{xNrN}, we have
\begin{eqnarray}
\label{max} \max_{1\le n \le N+1}\ln^+|(x_n^{-1})_{i,j}|_p & \le &
2\sum_{n=1}^N(|\ln |(g_{n+1}^{-1})_{i,i}|_p|+|\ln c_n|) \\
                & + & \nonumber \max_{i+1\le k\le j} \max_{1\le n \le
                N}\ln^+|(x_n^{-1})_{k,j}|_p.
\end{eqnarray}
Now by an elementary induction on $(j-k)\in [0,\dots,j-i]$ (with
$j$ fixed), we get from \eqref{max}
\begin{eqnarray}
\label{domination} \forall N\ge 1\quad
\frac{1}{N}\E\Big[\max_{1\le n \le N+1}
\ln^+|(x_n^{-1})_{i,j}|_p\Big]\le C K_p, \end{eqnarray} for some
constant $C>0$. Let now
$$\alpha_n^p:= \frac{1}{n} \ln^+|(x_n^{-1})_{i,j}|_p.$$
Again from \eqref{max} we get by induction on $j-k$ that a.s. for
all $p\notin P$, $(\alpha_n^p - CK_p)^+$ tends to $0$, when $n\to
+\infty$. By \eqref{domination} and Lebesgue theorem, we have even
that $\sum_p (\alpha_n^p -CK_p)^+$ converges to $0$ in $L^1$. So
for some constant $C'>0$,
\begin{eqnarray}
\label{partiehorsP} \pp\Big[\frac{1}{n}\sum_{p\notin P}
\lan x_n^{-1}\delta_n\ran^+_p \le C'\sum_{p\notin P}K_p \Big]\to 1.
\end{eqnarray}
\textit{Second step: the sum over $p\in P$.}
\newline
We will show now that for all $i\le j$ and all $p\in P$,
$$\left|\frac{(x_n^{-1}Z^p)_{i,j}}{(x_n^{-1})_{j,j}}\right|_p\le
e^{n\epsilon}$$ for $n$ large enough. Together with
\eqref{partiehorsP} this will conclude the proof of the
proposition. Without loss of generality we can always suppose that
$i=1$ and $j=d$. For $n\ge 1$, and $l\ge 1$, we have
$$\frac{(x^{-1}_nZ^p)_{l,d}}{(x^{-1}_n)_{d,d}}=\sum_{k=l}^d
\frac{(x^{-1}_n)_{l,k}}{(x^{-1}_n)_{d,d}}Z^p_{k,d}:=u_n^l$$ Using
that $x_nx^{-1}_n=\textrm{Id}$, we get by an elementary induction
on $l$ that, for $l\ge 1$,
\begin{eqnarray}
\label{unl}
u_n^l=\frac{(x_n)_{d,d}}{(x_n)_{l,l}}Z^p_{l,d}-\sum_{k=l+1}^d\frac{(x_n)_{l,k}}{(x_n)_{l,l}}u_n^k.
\end{eqnarray}
Now for $l<k$ let $A(l,k):=((x_n)_{i,j})_{l\le i<k,l<j\le k}$.
Next we will need the elementary formula
\begin{eqnarray} \label{determinantA} \frac{\det
A(l,k)}{\prod_{l'=l+1}^{k}(x_n)_{l',l'}}=\sum_{l'=l+1}^{k}(-1)^{l'-l+1}\frac{(x_n)_{l,l'}}{(x_n)_{l',l'}}\frac{\det
A(l',k)}{\prod_{l''=l'+1}^{k}(x_n)_{l'',l''}},
\end{eqnarray}
where by convention $A(k,k)=(1)$. We denote also for $l<k_1 < k$,
by $A(l,\widehat{k_1},k)$ the matrix $A(l,k)$, where the
$k_1^{th}$ line and the $(k_1-1)^{th}$ column are omitted. With
evident notation we define analogously
$A(l,\widehat{k_1},\dots,\widehat{k_r},k)$ for $l<k_1<\cdots
k_r<k$. For any $l<d$ we set $I_d^l:=\{l\}\cup(J-\{d\})$,
$J_d^l:=\{j\in J/\ j\neq l\}$ and
$$S_n^l:=\frac{\epsilon_l^p}{\prod_{j\in J}(x_n)_{j,j}}\det \Big( ((x_n)_{i,j})_{(i,j)\in I_d^l\times
J_d^l}\Big),$$ where $\epsilon_l^p=(-1)^{\text{Card}\{l<i<d\mid
\phi_p(i)\ge\phi_p(d)\}}$. By convention we set also $S_n^d=-1$.
We will need the
\begin{lem}
\label{zk} When $l\notin J$,
$$Z^p_{l,d}=\lim_{n\to \infty} S_n^l.$$
\end{lem}
We postpone the proof of this lemma to the appendix. Let
$\{i_1,\dots,i_s\}=J^c$ and for $l\ge 1$, let $k_l=\min \{k\le
s\mid i_k\ge l\}$. First we prove by induction on $d-l\ge 0$, that
\begin{eqnarray}
\label{recfacile} u_n^l&=&
-\mathbf{1}_{(l\in J)}\frac{(x_n)_{d,d}}{(x_n)_{l,l}}S_n^l\\
&+& \nonumber
\sum_{k=k_l}^s(-1)^{i_k-l}\frac{(x_n)_{d,d}^{i_k-l+1}\det
A(l,i_k)}{(x_n)_{l,l}\cdots(x_n)_{i_k,i_k}}(Z^p_{{i_k},d}-S_n^{i_k}),
\end{eqnarray}
where we recall our convention $A(l,i_{k_l})=(1)$ if $i_{k_l}=l$
(i.e. if $l\notin J$). In fact the result is trivial for $l=d$.
Now we suppose that it is true for $l$ strictly greater than some
$l_0$. Then Formula (\ref{recfacile}) for $l_0$ is a direct
consequence of (\ref{unl}) and (\ref{determinantA}), which proves
the induction step. Next we prove also by induction on $d-l\ge 0$,
that
\begin{eqnarray}
\label{formulerec} u_n^l & = &
\frac{(x_n)_{d,d}}{(x_n)_{l,l}}(Z^p_{l,d}-S_n^l)\\
&+& \nonumber \sum_{k=k_l}^s \frac{(-1)^{i_k-l}\det
A(l,\widehat{i_{k_l}},\dots,\widehat{i_{k-1}},i_k)}{(x_n)_{l,l}
\cdots(\widehat{x_n})_{i_{k_l},i_{k_l}}\cdots(\widehat{x_n})_{i_{k-1},i_{k-1}}\cdots(x_n)_{i_k-1,i_k-1}}u_n^{i_k},
\end{eqnarray}
where the notation $\hat{x}$ means that $x$ is omitted in the
list. Formula (\ref{formulerec}) is true for $l=d$. So we suppose
that it is true for $l$ strictly greater than some $l_0$. Then
observe that for any $l < k'<k$,
$$\det A(l,k)=\det A(l,k')\det A(k',k) +
(-1)^{k'-l}(x_n)_{k',k'}\det A(l,\widehat{k'},k).$$ Injecting this
in (\ref{recfacile}) and using the induction hypothesis we get
(\ref{formulerec}) for $l_0$, and we can conclude by the induction
principle. Eventually we prove again by induction on $d-l\ge 0$,
that $|u_n^l|_p\le e^{n\epsilon}$ for $n$ large enough. We suppose
that it is true for $l$ strictly greater than some $l_0$. For any
$l$ and any $k>k_l$, $\det A(l,\dots,\widehat{i_{k-1}},i_k)$ is
equal to the component on $e_l\dots
\wedge\widehat{e_{i_{k-1}}}\dots\wedge e_{i_k-1}$ of
$(ae_2\dots\wedge\widehat{ae_{i_{k-1}}}\dots\wedge ae_{i_k})$.
Therefore as in the proof of Lemma \ref{convergence}, we see that
$$\Big|\frac{\det
A(l,\widehat{i_{k_l}},\dots,\widehat{i_{k-1}},i_k)}{(x_n)_{l,l}
\cdots(\widehat{x_n})_{i_{k_l},i_{k_l}}\cdots(\widehat{x_n})_{i_{k-1},i_{k-1}}\cdots(x_n)_{i_k-1,i_k-1}}\Big|_p\le
e^{n\epsilon},$$ for n large enough. Moreover if $l_0\in J$, in
which case $Z^p_{l_0,d}=0$, we have also by the same argument
$|\frac{(x_n)_{d,d}}{(x_n)_{l_0,l_0}}S_n^{l_0}|_p\le
e^{n\epsilon}$ for $n$ large enough. Then we immediately prove the
result for $l_0$, by using the induction hypothesis and Formula
(\ref{formulerec}). This finishes the proof of the lemma.
\end{proof}

\noindent We are now ready for the

\vspace{0.2cm} \noindent \textbf{Proof of Proposition \ref{prop3}:}
Let $P$ be some finite subset of $\kP^*$ containing
$\infty$. For $z\in \prod_{p\in\kP^*} B_p$, let $z^P$ be its
natural projection on $\prod_{p\in P} B_p$. Let $K= \epsilon+C\sum_{p\notin
P} K_p$, where $\epsilon$, $C$ and $K_p$ are as in Lemma
\ref{estimgauge}. Then by Lemma \ref{estimgauge} (remember also \eqref{defgauge})
\begin{eqnarray}
\label{convl1} \pp\left[ x_n\in
\kG_{nK}^{\pi_n^P(Z^P_\infty)}\right]=\int_{\mathbf{B}}
\pp_n^z\left[\kG_{nK}^{\pi_n^P(z^P)}\right]d\nu(z) \to 1.
\end{eqnarray} Remember that $h^z$ denotes the
$\pp^z$-almost sure limit of $-\ln \pp_n^z(x_n)/n$. Consider
the set
$$A_n=\{g\in A(\Q)\mid -h^z-\epsilon <\ln\pp_n^z(g)/n<-h^z+\epsilon \}.$$
Then
$$\pp_n^z(A_n\cap \kG_{nK}^{\pi_n^P(z^P)})\le
e^{n(\epsilon-h)}\text{Card} \left(\kG_{nK}^{\pi_n(z^P)}\right)
\le e^{n(\epsilon-h^z)}e^{C'nK},$$ where $C'$ is the parameter of
the exponential growth of the gauges $(\kG^g)_{g\in H}$. Thus we
must have $C'K-h^z+\epsilon \ge 0$ for $\nu$-almost all $z\in \mathbf{B}$. Otherwise this would contradict
\eqref{convl1}. Since $\epsilon$ was arbitrarily chosen, we get
$$h^z\le C'C\sum_{p\notin P}K_p.$$ Letting now $P$ grow to
$\kP$, we obtain $h^z=0$, which concludes the proof of the proposition.
\hfill $\square$

\section{The case of a number field}
\label{secnf}
In this section $\K$ denotes a number field, i.e. a finite
extension of $\Q$. We refer to \cite{Sam} \cite{Ser} \cite{Wei}
for the general theory. Let $\kO$ be the ring of integers of $\K$.
The main difference with the rational case is that except for $\Q$
or imaginary quadratic extensions of $\Q$, the set $\kO^*$ of
units (the invertible elements of $\kO$) of $\K$ is infinite. So
we have to be careful when defining the gauges, to keep them with
uniform exponential growth. Namely we have to define $\lan k\ran$ for
$k\in \K$, in a suitable way. More precisely, let $\kP$ be the set
of prime ideals of $\kO$, and for $\kp \in \kP$ let $v_\kp$ be the
associated discrete valuation. Let $N_\kp=\text{Card}(\kO/\kp)$.
Following a usual convention (see e.g. \cite{Lan}), we define the
norm associated to $\kp$ by $$|k|_\kp:=N_\kp^{-v_\kp(k)},$$ for
all $k\in \K^*$. Let $N$ be the norm function on $\kO$. If
$k=x^{-1}y$ with $x,y\in \kO$ such that $(x)\wedge (y)=1$, then
$$\sum_{\kp\in \kP}|\ln |k|_\kp|=\ln|N(x)|+\ln|N(y)|.$$ Remember
that the norm of any unit is equal to $\pm 1$. Thus with the
previous notation we can not define $\lan k\ran$ as the sum
$\ln|N(x)|+\ln|N(y)|$ like in the rational case. Otherwise the
associated gauges would have an infinite cardinality. So we have to add a term corresponding to archimedean norms. Remember that $V_\infty$ denotes the set of norms extending
the usual absolute value $|\cdot|$ on $\Q$. If $v\in V_\infty$ we
will write (with a slight abuse of notation) $v"="|\cdot|_v$ and
we define the norm $|\cdot|_v:=|\cdot|_v^{\epsilon_v}$ on $\K$,
where $\epsilon_v=1$ if $\K_v=\R$ whereas $\epsilon_v=2$ if
$\K_v=\C$. Then we have the product formula (see
\cite{Lan})
$$\prod_{\kp\in \kP}|k|_\kp \times \prod_{v\in V_\infty}|k|_v = 1,$$
for all $k\in K^*$, which implies by the way the identity
\begin{eqnarray}
\label{identitephi} \sum_{\kp\in \kP}\phi_\kp+\sum_{v\in
V_\infty}\phi_v=0.
\end{eqnarray}
Now we fix some archimedean norm $|\cdot|_{v_0}$ and we define
$$\lan k\ran :=\sum_{\kp\in \kP}|\ln |k|_\kp|+\sum_{v\neq v_0} |\ln
|k|_v|.$$ In this way the set of $k\in \K^*$ such that $\lan k\ran \le
C$ has a cardinality bounded by $\text{const}\cdot
e^{\text{const}\cdot C}$, where the constants are independent of
$C$. Then we can define the height function on the adele ring and
the associated gauges, in the same way as in the rational case.
Now the only other change in the proof is the definition of the
$q_n^i$ (see section \ref{secestimate}). Remember that the set of
units is isomorphic to $\Z^{r_1+r_2-1}\times G$, where $G$ is
cyclic, $r_1$ is the number of embedding of $\K$ in $\R$ and
$2r_2$ the number of embedding in $\C$. We set
$$q_n^i=\prod_{\kp\in\kP}p^{-[\frac{n\phi_{\kp}(i)}{\ln
p}]}\prod_{v\neq v_0}u_v^{-[n\alpha_v]},$$ where in the first
product, for each $\kp$ the prime number $p$ is such that $v_\kp$ extends
$v_p$ on $\Q^*$, and in the second product the $u_v\in \kO^*$ and
the $\alpha_v\in \R$ are chosen as follows. For $(u_v)_{v\neq
v_0}$ take any basis of $\Z^{r_1+r_2-1}$ (seen as a subset of
$\kO^*$). Then the matrix $(\ln ||u_v||_w)_{v,w\neq v_0}$ is
invertible (see the proof of Theorem $1$ p.$72$ in \cite{Sam}). So
one can choose $(\alpha_v)_{v\neq v_0}$ such that
$$\sum_{v\neq v_0} \alpha_v \ln ||u_v||_w=\sum_{v\neq w} \phi_v(i),$$
for all $w\neq v_0$. Thus with \eqref{identitephi} one can check
that the analogue of Lemma \ref{qni} holds. The other parts of the
proof are unchanged. We leave the details to the reader. 

\section{Appendix}
\noindent \textbf{Proof of Lemma \ref{zk}:} Let $l\notin J$. Let
$i_1<\cdots<i_r$ be the elements of $I_d^l$. Assume that
$e_{i_1}\wedge \cdots e_{i_r}$ is the $k^{th}$ element of the
basis $\kB_r$ (with the notation of section \ref{secmuboundary}).
First by definition of $Z^p$, we can see that
$Z^p_{l,d}=\epsilon^p_l Z_k^p(d)$. Next we have seen in the proof
of Lemma \ref{convergence} that $Z^p_k(d)=\lim_{n\to \infty}
(x'_n)_{k,m}$. We will show in fact directly that for any $a\in
A(\Q)$,
$$a^{(r)}_{k,m}=\det\Big(\underbrace{ (a_{i,j})_{(i,j)\in I_d^l\times
J}}_{:=M(l)}\Big).$$ Naturally it will imply the lemma. We prove
the result by induction on $h=d-l$. If $h=1$, i.e. $l=d-1$, then
$I_d^l=\{j_1,\dots,j_{r-1},d-1\}$, and
$a^{(r)}_{k,m}=a_{j_1,j_1}\dots a_{j_{r-1},j_{r-1}}a_{d-1,d}$.
Then the result is immediate. We prove now the induction step from
$h$ to $h+1$. We suppose that $j_{s-1}<l<j_s$ for some $s$ (if
$s=1$ we have just $l<j_s$). The coefficient $a^{(r)}_{k,m}$ is
equal to the component on $e_{j_1}\dots \wedge e_l\dots\wedge
e_{j_{r-1}}$ of $(ae_{j_1}\wedge \dots \wedge ae_{j_r})$. This
component is equal to the sum over $k\in [s,\dots,r]$, of the
components of $(\prod_{j<l}a_{j,j})(e_{j_1}\dots \wedge
ae_{j_s}\dots \wedge a_{l,j_k}e_l\dots
\wedge\widehat{ae_{j_k}}\dots \wedge ae_{j_r})$ on $e_{j_1}\dots
\wedge e_{j_k}\dots \wedge e_{j_{r-1}}$. But by the induction
hypothesis, for each $k\in [s,\dots,r]$, the corresponding
component is equal to $a_{l,j_k}$ times the cofactor of
$a_{l,j_k}$ in the matrix $M(l)$. This gives exactly the formula
of the determinant of $M(l)$. Therefore the proof of the lemma is
finished. \hfill $\square$

\vspace{0.5cm}

\noindent Universit\'e d'Orl\'eans, F\'ed\'eration Denis Poisson, Laboratoire MAPMO \\
B.P. 6759, 45067 Orl\'eans cedex 2, France.

\vspace{0.2cm}

\noindent current address:  D\'epartement de Math\'ematiques, B\^at. 425, Universit\'e Paris-Sud, F-91405 Orsay
cedex, France. 
e-mail: bruno.schapira@math.u-psud.fr


\begin{thebibliography}{99}
\bibitem{BS} \textbf{Bader U., Shalom Y.:} \textit{Factor and normal subgroup theorems for lattices in products of groups},
Invent. Math. 163, (2006), 415--454.
\bibitem{Bro} \textbf{Brofferio S.:} \textit{The Poisson Boundary of
random rational affinities},  Ann. Inst. Fourier 56, (2006),
499--515.
\bibitem{CKW} \textbf{Cartwright D. I., Kaimanovich V. A., Woess
W.:} \textit{Random walks on the affine group of local fields and
of homogeneous trees}, Ann. Inst. Fourier 44 (1994), 1243--1288.
\bibitem{Der} \textbf{Derriennic Y.:} \textit{Entropie, th\'eor\`emes
limite et marches al\'eatoires}, in Probability measures on groups
VIII (Oberwolfach, 1985), LNM 1210, pp. 241--284, Springer,
Berlin, (1986).
\bibitem{Eli} \textbf{Elie L.:} \textit{Noyaux potentiels associ\'es aux
marches al\'eatoires sur les espaces homog\`enes. Quelques exemples
clefs dont le groupe affine}, in Th\'eorie du potentiel (Orsay,
1983), volume 1096 of Lectures Notes in Math., 223--260, Springer,
Berlin, 1984.
\bibitem{Fur} \textbf{Furman A.:} \textit{Random walks on groups
and random transformations}, Handbook of dynamical systems, vol.
1A, pp. 931--1014, Amsterdam: North-Holland (2002).
\bibitem{F} \textbf{Furstenberg H.:} \textit{A Poisson formula for
semi-simple Lie groups}, Ann. of Math. 77 (1963), 335--386.
\bibitem{F2} \textbf{Furstenberg H.:} \textit{Boundary theory and
stochastic processes on homogeneous spaces}, in Harmonic analysis
on homogeneous spaces (Proc. Sympos. Pure Math., Vol. XXVI,
Williams Coll., Williamstown, Mass., 1972), p. 193--229, Amer.
Math. Soc., Providence, R.I., (1973).
\bibitem{GR} \textbf{Guivarc'h Y., Raugi A.:} \textit{Fronti\`ere de Furstenberg, propr\'et\'es de contraction et th\'or\`emes de convergence},
(French) Z. Wahrsch. Verw. Gebiete 69, (1985), 187--242.
\bibitem{Kai} \textbf{Kaimanovich V. A.:} \textit{The Poisson formula for groups with hyperbolic properties}, Ann. of Math. (2) 152, (2000), 659--692.
\bibitem{KV} \textbf{Kaimanovich V. A., Vershik  A. M.:}
\textit{Random walks on discrete groups: boundary and entropy},
Ann. Probab. 11, (1983), 457--490.
\bibitem{Lan} \textbf{Lang S.:} \textit{Introduction to diophantine approximations}, Addison-Wesley Publishing Co.,
Reading, Mass.-London-Don Mills, Ont., (1966).
\bibitem{Led} \textbf{Ledrappier F.:} \textit{Poisson boundaries of discrete groups of matrices},
Israel J. Math. 50, (1985), 319--336.
\bibitem{Mar} \textbf{Margulis G. A.:} \textit{Discrete subgroups of semisimple Lie groups},
Ergebnisse der Mathematik und ihrer Grenzgebiete 17.
Springer-Verlag, Berlin, (1991), x+388 pp.
\bibitem{Rau} \textbf{Raugi A.:} \textit{Fonctions harmoniques sur les groupes localement compacts \`a base d\'enombrable},
Bull. Soc. Math. France M\'em. No. 54, (1977), 5--118.
\bibitem{Sam} \textbf{Samuel P.:} \textit{Th\'eorie alg\'ebrique des nombres (French)}, Hermann, Paris (1967), 130
pp.
\bibitem{Ser} \textbf{Serre J.P.:} \textit{Corps locaux (French)}, Sec. edition, Publications  of University Nancago, No. VIII. Hermann, Paris, (1968), 245 pp.
\bibitem{War} \textbf{Warner G.:} \textit{Harmonic analysis on semi-simple Lie groups. I}, Die Grundlehren der mathematischen Wissenschaften, Band 188.
Springer-Verlag, New York-Heidelberg, (1972), xvi+529 pp.
\bibitem{Wei} \textbf{Weil A.:} \textit{Basic number theory}, Third edition, Die Grundlehren der Mathematischen Wissenschaften, Band 144.
Springer-Verlag, New York-Berlin, (1974), xviii+325 pp.

\end{thebibliography}
\end{document}